\mag 1200
\input amstex
\input amsppt.sty

\overfullrule 0pt

\hsize 6.25truein
\vsize 9.63truein

\let\<\langle \let\>\rangle

 \let\eps\varepsilon \let\epsilon\eps

 \let\phi\varphi

\def\C{\Bbb C}
\def\R{\Bbb R}
\def\Z{\Bbb Z}

\def\F{\Cal F}

  \def\Fwh{\,\wh{\!\F}}

\def\A{\roman A}

\def\E{\roman E}

\let\Fun\Fwh

% my notations

\def\A{\Cal A}

\def\Fun/{ {\text{Fun}}}
\def\End/{ {\text{End}}}
\def\Hom/{ {\text{Hom}}}
\def\Rm/{\^{$R$-}matrix}

\def\M/{\Cal M }
\def\A/{$A_{\tau,\eta}(sl_2)$}
\def\Am/{$A_{\tau,\eta,\mu}(sl_2)$}

\def\E/{$E_{\tau,\eta}(sl_2)$}

\def\Ex/{ \widetilde {  E   }}
\def\EE/{ $\widetilde {  E   }_{\tau,\eta}(sl_2)$}

\def\eqg/{elliptic quantum group}
\def\Em/{$E_{\tau,\eta,\mu}(sl_2)$}
\def\E/{$E_{\tau,\eta}(sl_2)$}

\def\Ce/{\C/{1\over \eta}\Z }

\topmatter

\title
Quantization of geometry associated to the quantized Knizhnik-Zamolodchikov
equations $^*$
\endtitle
                                
\author
Alexander Varchenko
\endauthor

\thanks
The author was supported in part by NSF grant
DMS-9501290 \newline $\,^*$ This is a talk given at the
1995  Les Houches School on Quantum Symmetries.
\endthanks

\date
May, 1996
\enddate

\abstract

It is known that solutions of the  Knizhnik-Zamolodchikov 
differential 
equations
are given by integrals of closed 
differential forms over suitable cycles. In this paper a quantization
of this geometric construction is described leading to solution of the 
quantized Knizhnik-Zamolodchikov difference equations.

\endabstract

 \leftheadtext{ A.Varchenko }
  \rightheadtext{ Quantization of geometry associated to qKZ}
\endtopmatter

\document

\head 1. Introduction \endhead

The Knizhnik-Zamolodchikov equations (KZ) are fundamental differential
equations discovered in conformal field theory by Knizhnik and Zamolodchikov
in the beginning of the 80-s. It is known now  that the KZ equations
are  the same equations 
as the differential equations for multidimensional
hypergeometric functions. Hypergeometric functions are purely geometric
objects and there is a certain geometry associated to these functions.

In mathematical physics the KZ differential equations were quantized.
The quantized Knizhnik-Zamolodchikov equations (qKZ) are important difference 
equations [FR], [JM], [S]. 
It turns out that there is a geometric theory
of the qKZ equations which is  a quatization of the geometric
theory of hypergeometric functions.
This quantization of geometry is the subject of my talk.

First I will
 sketch the geometry of hypergeometric functions and then its quantization.

In this lecture I
will tell about my joint work with V.Tarasov [TV2].

The author thanks the organizers of the school, and
in particular K.Gawedzki, for invitation to give
this talk.

\head 2 KZ equations \endhead

Let $ \frak g$ be a simple Lie algebra, for instance $sl_2$.
Let
$V_1,\ldots,V_n$ be $\frak g$-modules, $V=V_1\otimes \ldots \otimes V_n$. 
Let $\Omega \in \frak g \otimes \frak g$ be the tensor corresponding to
an invariant scalar product. For $ i<j$, let $\Omega_{i,j}:V\to V$ be
the linear operator acting  as $\Omega$
 in $V_i\otimes V_j$ and as the identity operator in all of the
other factors.
 The KZ equations on 
an $V$-valued function $\Psi (z_1,\ldots,z_n)$ have the form
$$
{\partial \Psi \over \partial z_i} = {1\over \kappa}\,
\sum_{j,j\neq i}{\Omega_{i,j}\over z_i-z_j}\, \Psi, \qquad i=1,\ldots,n.
$$
Here $\kappa$ is a parameter of the equation. 

The KZ differential equations have very rich mathematical structures,
they are closely connected to affine Lie algebras,  quantum groups,
topology of knots and threefolds.

\head 3. Hypergeometric functions \endhead

There is a geometric source of differential equations of this type.
These are differential equations for multidimensional hypergeometric 
functions.

\proclaim {Claim [SV, V]}

The differential equations arising as differential equations 
for hypergeometric functions
are "precisely" the same as the KZ equations.
\endproclaim

{\it Example\/} of differential equations arising
as differential equations for hypergeometric functions.

Fix numbers $a_1,\ldots,a_n, \kappa$. Let
$$
\Phi(t,z_1,\ldots,z_n) = \prod^n_{\ell=1}
(t-z_\ell)^{a_\ell/\kappa}.
$$
For fixed $z_1,\ldots,z_n$,
 consider the complex line $\C$ and denote by $\gamma_m$
the oriented interval in $\C$ from $z_m$ to $z_{m+1}$, $m=1,\ldots,n-1$.

Let
$$
\Psi(z)=\Psi_\gamma(z_1,\ldots,z_n)=
(\int_\gamma \Phi\,{dt\over t-z_1},\ldots,\int_\gamma \Phi\,{dt\over t-z_n})
$$
where $\gamma=\gamma_1,\ldots,\gamma_{n-1}$.
Then
$$
{\partial
\Psi \over \partial z_i} = 
{1\over \kappa}\,\sum_{j,j\neq i}{\Omega_{i,j}\over z_i-z_j}\, \Psi, \qquad i=1,\ldots,n,
$$
where $\Omega_{i,j}(a)$ are some matrices independent on $z$ and $\gamma$.
Different $\gamma$ give different solutions to the same equations.

\head Geometry of hypergeometric functions \endhead

For fixed $z_1,\ldots,z_n$, consider a complex
$$
0\, \to\, \Omega^0\, \to\, \Omega^1\, \to \, 0,
$$
where 
$$
\Omega^0=\{f(t)\Phi(t,z)\,|\, f \text{ is any rational function
regular on}\ \, \C-\{z_1,\ldots,z_n\}\},
$$
$$
\Omega^1=\{f(t)\Phi(t,z)dt\,|\, f \text{ is any rational function
regular on}\ \,\C-\{z_1,\ldots,z_n\}\},
$$
and the differential of the complex is the standard differential.

\proclaim{Theorem}

For generic $a_1,\ldots,a_n,\kappa$, we have $H^0 = 0,$
dim $H^1\,=\, n-1$, and the differential forms
$\Phi dt/(t-z_1),\ldots,\Phi dt/(t-z_{n-1})$ form a basis in $H^1$.
 
\endproclaim

Let $H_1=(H^1)^*$ be the dual space, the first
homology group. Each interval $\gamma_m, \,
 m=1,\ldots,n-1$, defines a linear function on $\Omega^1$,
$$
[\gamma_m] \,:\, f\Phi dt \mapsto \int_{\gamma_m} f\Phi dt.
$$
By the Stokes' theorem 
we have $\int_{\gamma_m}d(f\Phi) = 0$ and , hence, 
$[\gamma_m]$ defines an element
of $H_1$.

\proclaim{Theorem}

The elements $[\gamma_m], m=1,\ldots.n-1$, form a basis in $H_1$.

\endproclaim

This fact is based on the formula
$$
{\text {det}}_{1\leq \ell, m \leq n-1 }
({a_\ell \over \kappa} \int_{\gamma_m} \Phi \, {dt\over
t-z_\ell }) =  {\Gamma({a_1\over \kappa}+1)\ldots
\Gamma({a_n\over \kappa}+1)\over \Gamma({a_1+\ldots+a_n\over \kappa}+1)}
\prod_{i\neq j}\,(z_i-z_j)^{a_j/\kappa}.
$$
For $n=2,\, (z_1,z_2) =(0,1)$, the formula takes the form
$$
b\int^1_0 t^a(1-t)^{b-1}dt = {\Gamma(a+1)\Gamma(b+1)\over
\Gamma(a+b+1)}.
$$

We have a correspondence: conformal field theory corresponds 
to geometry 
of hypergeometric functions,
the space of conformal blocks corresponds to the cohomology space
$H^1$ associated to the function 
$\Phi(t,z)$, the KZ equations correspond to the differential
equations
for integrals of basic closed differential forms
(to the Gauss-Manin connection), solutions of the KZ equations correspond to
cycles $\gamma \in H_1$.

\subhead Remark \endsubhead
More general solutions of the KZ equations are constructed
geometrically in a similar way starting from a more general function
$$
\Phi (t_1,\ldots,t_k,z_1,\ldots,z_n) \, = \, \prod (t_i-z_m)^{a_{i,m}/
\kappa} \,\prod (t_i-t_j)^{b_{i,j}/\kappa}.
$$

\head 5. qKZ equations\endhead

The KZ differential equations were quantized in  
[S], [FR], [JM]. The qKZ equations are difference equations satisfied by
form factors in [S], by matrix elements of intertwining operators in
[FR], by correlation functions of statistical  models
in [JM].

There are three versions of the qKZ equations: rational, trigonometric 
and elliptic.
We discuss here the rational version.
Let $V_1,\ldots,V_n$ be $sl_2$ modules with
highest weights, $V=V_1 \otimes \ldots
\otimes V_n$.
The rational qKZ on an $V$-valued function $\Theta(z_1,\ldots,z_n)$
 has the form
$$
\Theta (z_1,\ldots,z_m+p,\ldots,z_n) = K_m(z_1,\ldots,z_n,p)\,\Theta
(z_1,\ldots,z_n), \qquad m=1,\ldots,n,
$$
where $p$ is the step of the equations, $K_m:V\to V$ are linear 
operators defined as follows.

For any $i,j$ there is a a matrix called the
rational $R$-matrix
$$
R_{V_i,V_j}(x)\, :\, V_i\otimes V_j \, \to \, V_i \otimes V_j
$$
sutisfying the Yang-Baxter equation
$$
R_{V_i,V_j}(x)R_{V_i,V_k}(x+y)R_{V_j,V_k}(y)=
R_{V_j,V_k}(y)R_{V_i,V_k}(x+y)R_{V_i,V_j}(x)
$$
and normalized by the condition
$$
R_{V_i,V_j}(x) \, :\, v_i \otimes v_j\, \to\,
v_i \otimes v_j\, ,
$$
where $v_i, v_j$ are highest weight vectors,
see [JKMO],[KRS]. Then
$$
K_m = R_{V_m,V_{m-1}}(z_m-z_{m-1}+p) \ldots R_{V_m,V_1}(z_m-z_1+p)\,
R_{V_m,V_n}(z_m-z_n)\ldots R_{V_m,V_{m+1}}(z_m-z_{m+1}).
$$

The rational qKZ difference equation turns into the 
KZ differential equation under the following limiting procedure.

Let $z_m=SZ_m$, where $S\gg 1$ is a number, $Z_m$ are new variables.
Then for the new variables and a new function $\Psi(Z):=\Theta
(Sz)$ we have $\Psi(\ldots,Z_m +p/S,\ldots)
=K_m(SZ,p)\,\Psi(Z_1,\ldots,Z_n)$ and
$$
K_m(SZ,p)=1 + {1\over S}\,\sum_{j\neq m}\,{\Omega_{m,j}+c_{m,j}
\over
Z_m-Z_j}
+\Cal O({1\over S^2}),
$$
where $c_{m,j}$ are suitable numbers.
Hence, the difference equation turns into the differential equation
$$
{\partial \Psi\over 
\partial Z_m} = {1\over p}\,\sum_{j\neq m}\,{\Omega_{m,j}+
c_{m,j}\over
Z_m-Z_j}\, \Psi
$$
when $S$ tends to infinity.

\head 6. Solutions to the QKZ equations and eigenvectors of commuting
Hamiltonians \endhead

Assume that $p$ tends to zero and
$$
\Theta(z,p)=e^{S(z)/p}(f(z) + p f_1(z) + \ldots).
$$
Then
$$
K_m(z,p=0) f(z) =\lambda_m(z) f(z),
$$
where $\Lambda_m(z)=e^{\partial S(z)/\partial z_\ell}$.
This means that
the leading term $f(z)$ is an eigenvector of the commuting
operators $K_m(z,p=0) , m=1,\ldots,n$. The operators
$K_m$ are Hamiltonians of a quantum spin chain model .
Having a solution to the qKZ equations and computing its
quasiclassical asymptotics one can construct eigenvectors
of Hamiltonians. These eigenvectors coincide with the Bethe vectors
constructed by the Bethe ansatz [TV1].

\head 7. Solutions to the  qKZ equations\endhead

As we know the KZ differential equations are realized as differential equations
for closed differential forms over cycles depending on parameters.
It turns out that there is a quantization of all these
 geometric objects leading to solutions of the
qKZ difference equations.

Let us describe a
 $p$-analog of the function
$\Phi(t,z)=\prod(t-z_\ell)^{a_\ell/\kappa}$

The function $T^{a/\kappa}$ satisfies the differential equation
$$ 
{dY\over dT}={a\over \kappa} \,Y.
\tag 1
$$
A $p$-analog of this differential equation is
the difference equation
$$ 
y(t+p)={t+a\over t-a} \,y(t).
\tag 2
$$
If $t=ST,\, Y(T):=y(ST)$, then the difference equation takes the form
$$ 
Y(T+{p\over S})=(1+{2a\over S}{1\over T -a/S}) \,Y(T)
$$
and turns into the differential equation (1) with $\kappa =p/2$ as $S$ tends
 to infinity.

Equation (2) has a solution
$$
y=\Gamma({t+a\over p})\Gamma(1-{t-a\over p})
e^{-\pi i t/p}.
$$
Introduce a function
$$
\Phi_p(t,z)=\prod^n_{\ell=1} \Gamma({t-z_\ell+a_\ell\over p})
\Gamma(1-{t-z_\ell-a_\ell\over p})
e^{-\pi i (t-z_\ell)/p}.
$$
This function is a $p$-analog of the function $\Phi$ of section 3.

\subhead Properties \endsubhead

1.
$$
\Phi_p(t+p,z) = \prod_{\ell=1}^n{t-z_\ell+a_\ell \over
t-z_\ell-a_\ell } \,
\Phi_p(t,z)=:b_0(t,z)  \Phi_p(t,z).
$$

2.
$$
\Phi_p(t,z_1,\ldots.,z_\ell + p,\ldots,z_n)=
{t-z_\ell-a_\ell -p \over
t-z_\ell+a_\ell -p} \,
\Phi_p(t,z)=: b_\ell(t,z)   \Phi_p(t,z).
$$
We have
$$
b_m(\ldots,z_\ell +p,\ldots)\, b_\ell (z)=
b_\ell(\ldots,z_m +p,\ldots)\, b_m (z),
$$
where
$\ell, m \in \{0,\ldots,n\}$.

The singularities of $\Phi(t,z)$ for fixed $z$
are the points
$
\{t\in \C\, | \, t=z_\ell-a_\ell - Np, $ \newline $
t=z_\ell+a_\ell + (N+1)p, \,
{\text{
where}}\,  \ell = 1,\ldots,n, \, N = 0,1,2,\ldots\}. $

Denote by $Sing^\vee (z)$ the set $\{t\in \C\, | \,
t\in \C \, | \,
t=z_\ell-a_\ell + (N+1)p, $ \newline $
t=z_\ell+a_\ell - Np,\, {\text{where}} 
\, \ell = 1,\ldots,n, \, \ell = 0,1,2,\ldots\}. $

Denote by $F(z)$ the linear space of rational functions
$f:\C \to \C$ , regular on $\C- Sing^\vee(z)$ and
with at most
simple poles at
$Sing^\vee(z)$ .

Consider a complex 
$$ 
0\,\to F(z) \, \to \, F(z) \, \to 0, \,
f(t) \, \mapsto \, f(t+p) b_0(t,z) - f(t). 
$$

\subhead Property \endsubhead
 The differential of this complex, $D_p(z)$,
preserves the space $F(z)$. 

Introduce the $p$-cohomology group
by $H_{(p)}^1(z)=F(z)/D_p(z)F(z)$.

Introduce linear maps 
$$
B_\ell(z):F(\ldots,z_\ell+p,\ldots) \to F(\ldots,z_\ell,\ldots),\,
f(t) \mapsto f(t)b_\ell(t,z)
$$
where $\ell = 1,\ldots, n$.

\subhead Properties \endsubhead

1. $B_\ell(z)$ is an isomorphism.

2. $B_\ell(z)$ commutes with the differential $D_p.$

3. $B_m(\ldots,z_\ell+p,\ldots)B_\ell(z)=B_\ell(\ldots,z_m+p,\ldots)B_m(z). $ 

\proclaim {Corollary} There are well defined isomorphisms
$$
B_\ell(z) : H^1_{(p)}(\ldots,z_\ell+p,\ldots) \to
H^1_{(p)}(\ldots,z_\ell,\ldots).
$$
These isomorphisms satisfy property 3.
\endproclaim

\proclaim {Corollary} If $z-z' \in (p\Z)^n$, then
$H^1_{(p)}(z) $ and $H^1_{(p)}(z') $ are canonically identified.

\endproclaim

We  call these isomorphisms a {\it discrete flat connection}.

Let $H_1^{(p)}(z) $ be the space dual to $H^1_{(p)}(z) $.
The space will be called the {\it $p$-homology group},
 its elements will be called
 {\it $p$-cycles}.

Let $\gamma :z \mapsto \gamma(z) \in H_1^{(p)}(z) $
be a section. The section is called {\it periodic} if $B_\ell(z)\, 
\gamma(\ldots,z_\ell+p,\ldots)=\gamma(z)$. We define 
periodic sections with values in
$H^1_{(p)}(z)  $ analogously .

\proclaim {Theorem }

 For generic $z_1,\ldots,z_n,a_1,\ldots,a_n,p$,
the dimension of $H_1^{(p)}(z) $ equals $n-1$.
The elements
$$
w_j(z,t)={1\over t-z_j-a_j}
\prod_{\ell=1}^{j-1}
{t-z_\ell+a_\ell \over t-z_\ell-a_\ell } \,\in \, F(z),
\qquad j=1,\ldots,n-1,
$$
generate a basis in $H^1_{(p)}(z) $ . Under the above limit
$S \to \infty$ ,
the element $w_j$ tend to $1/(T-Z_j)$.
\endproclaim

\head 8. Difference equations of the discrete connection
\endhead

Since the elements $[w_1(z)]$,\,$\ldots$,\,$[w_{n-1}(z)]$
form a basis in $H_1^{(p)}(z) $ , the isomorphisms
$B_\ell(z)$ could be given by matrices $\beta_\ell(z)$,
$$
B_\ell(z)[w_j(\ldots,z_\ell+p,\ldots)]=w_i(z)\,\beta_{\ell,j}^i(z).
$$

For 
any $p$-cycle $\gamma \in H^1_{(p)}(z) $ 
introduce its coordinate vector
$$
I=(\<w_1(z),\gamma\>,\ldots,\<w_{n-1}(z),\gamma\>).
$$

\proclaim{Theorem}

Let a $p$-cycle $\gamma(z) \in H_1^{(p)}(z)$ be periodic,
then its coordinate
vector satisfies the system of difference
equations
$$
I(\ldots,z_\ell+p,\ldots)=I(z)\beta_\ell(z), \qquad, \ell =1,\ldots,n.
\tag 3
$$
Changing periodic cycle $\gamma$, one constructs all solutions to system (3).
\endproclaim

{\it Proof}. $I(\ldots,z_\ell +p,\ldots)=
(\ldots,\<w_j(\ldots,z_\ell+p,\ldots),\gamma
(\ldots,z_\ell+p,\ldots)\>,\ldots)=$\newline
$(\ldots,\<B_\ell(z) w_j(\ldots,z_\ell+p,\ldots),
(B_\ell(z)^*)^{-1} \gamma
(\ldots,z_\ell+p,\ldots)\>,\ldots)=$ \newline
$(\ldots,\<w_i(z)\,\beta_{\ell,j}^i,\gamma(z)\>,\ldots)=
I(z) \beta_\ell(z). $

\proclaim{Theorem}
System (3) coincides
with a special case of the qKZ equations. This special case
is the quantization of the special case of the KZ differential
equations associated to the function $\Phi(t,z)=
\prod (t-z_\ell)^{a_\ell/\kappa}$
of section 3.

\endproclaim

All qKZ equations associated to representations of $sl_n$ could be constructed 
in a similar way.

\subhead Remark \endsubhead
The theorem
says that the $p$-homology classes can be naturally identified 
with suitable vectors in the corresponding tensor product of
representations of $sl_2$ in such a way that the action 
of $R$-matrices is identified with the action of the discrete connection.

\head 9. $p$-Homology theory \endhead

According to our definitions
$$
D_p(z) : f(t) \mapsto (f(t+p)\Phi_p(t+p,z)-f(t)\Phi_p(t,z))/\Phi_p(t,z),
$$
$H^1_{(p)}(z)=F(z)/D_p(z)F(z)$ and
$H_1^{(p)}(z)=(H^1_{(p)}(z))^*$.
In this section we will construct linear
functions on $H^1_{(p)}(z)=F(z)/D_p(z)F(z)$ , that is, we will
construct elements of
$H_1^{(p)}(z)$.

\subhead A naive idea \endsubhead

For $\xi \in \C$ and a function $ h:\C \to \C$
with compact support define the Jackson integral by
$$
\int_{[\xi]_p} h\, d_pt = p
\sum^\infty_{\ell=-\infty} h(\xi +\ell p).
$$
\subhead Property \endsubhead Let $h$ be of the form
$h(t)=D_pg(t)=g(t+p)-g(t)$, then
$$
\int_{[\xi]_p} h\, d_pt = p
\sum^\infty_{\ell=-\infty} (g(\xi +(\ell+1) p)-g(\xi +\ell p))=0.
$$
Therefore, a naive way to construct elements of 
$H_1^{(p)}(z)$ would be to take  any $\xi$ and define for any
$f\in F(z)$
$$
<f,[\xi]_p> = \int_{[\xi]_p} \Phi_p(t,z)f(t) d_pt.
$$
Then such a linear function on $F(z)$ would be zero on 
$D_p(z) F(z)$ and would define an element of 
$H_1^{(p)}(z)$. 

Unfortunately, this idea does not work since
the integral does not converge. However, a certain modification
of the idea could be realized.

Choose $p$ so that $p$ is imaginary,
$p \in i\R$, and Im $p > 0$. Choose
$z_1,\ldots,z_n$ to be real, and choose
$ a_1,\ldots,a_n$ to be  imaginary
and lying in the upper half plane,
$a_\ell \in i\R$, and Im $a_\ell > 0$. 

Under these assumptions the real line $\R$
separates in $\C$ the poles of the factors
$\prod^n_{\ell=1} \Gamma({t-z_\ell+a_\ell\over p})$ and
the factors $\prod^n_{\ell=1} 
\Gamma(1-{t-z_\ell-a_\ell\over p})$ of the function $\Phi_p$.

Let $G_m:\C \to \C,\, t \mapsto e^{2\pi i mt/p},\, m=1,\ldots,n-1$.
The functions $G_m$ are the simplest $p$-periodic functions.

For every such $m$, consider a map
$$
[G_m]:F(z) \to \C, \, \, \, 
f \mapsto
\int_\R G_m(t) \Phi_p(t,z) f(t) dt.
$$
As we will see, 
the maps $[G_m]$ are $p$-analogs of
the intervals $\gamma_m$.

\subhead Properties of the linear functionals $[G_m]$
\endsubhead

1. The functionals $[G_m]$ , $m=1,\ldots,n-1$, are 
well defined. Moreover, for any other integer $m$,
the linear functional $[G_m]$ is not defined on $F(z)$.

{\it Proof.} The function $\Phi_p$ is a product 
of gamma functions. 
When $x\to\infty$, the gamma function $\Gamma (x)$ has asymptotics
$$
\Gamma (x) = x^{-1/2} e^{x( {\text{ln}} \,x\, -\,1)}(2\pi)^{1/2} + \ldots
$$ 
These asymptotics imply the property.

2. For $m=1,\ldots,n-1, \, [G_m]|_{D_p(z)F(z)} =0$.

{\it Proof}.
$$
\int_{\R} G_m(t)[\Phi_p(t+p,z)f(t+p)-\Phi_p(t,z)f(t)]\,dt =
$$
$$
\int_{\R+p} G_m(t)\Phi_p(t,z)f(t)\,dt -
\int_{\R} G_m(t)\Phi_p(t,z)f(t)\,dt = 0.
$$
Here we use periodicity of the function $G_m(t)$
and the fact that there are no poles of the integrand between
$\R$ and $\R+p$.

\subhead Corollary
\endsubhead The functional $[G_m]$ defines an element
of $H_1^{(p)}(z)$.

3. If $z_\ell = SZ_\ell$ and $S \to \infty$, then the
linear functional $[G_m]$
tends to the interval $[Z_m,Z_{m+1}]$
in the following sence:
$$
\int_{\R} G_m(t)\Phi_p(t,z)w_\ell(t)\,dt =
C_m(S,p)\,(\int_{Z_m}^{Z_{m+1}} \prod_{j=1}^n (T-Z_j)^{2a_j/p}
{dT\over T-Z_\ell} \,+\, O(1/S))
$$
for every $\ell$. Here $C_m$ is an explicitly given function.

{\it Proof\/} follows from the Stirling formula.

4.  The sections $[G_m](z) \in H_1^{(p)}(z) $ are $p$-periodic.

{\it Proof.} For $f \in F(\ldots,z_\ell +p,\ldots)$ we have
$\<[G_m](\ldots,z_\ell+p,\ldots), [f]\>:= $ \newline
$\int_{\R} G_m(t)\,\Phi_p(t,\ldots,z_\ell+p,\ldots)\,f(t)\,dt =
\int_{\R} G_m(t)\,\Phi_p(t,z)\,b_\ell(t,z)\,f(t)\,dt = $\newline
$\<[G_m](z),B_\ell(z)[f]\>.$

5. The elements $[G_m](z)$, $m=1,\ldots,n-1$, form a basis
in $H_1^{(p)}(z)$.

\subhead Corollary \endsubhead The vectors
$$
\Theta_m(z)=(\int_\R G_m \Phi_pw_1 dt, \ldots,
\int_\R G_m \Phi_pw_{n-1} dt), \qquad m=1,\ldots,n-1,
$$
form a basis of solutions to the difference equations of
the discrete connection on $H^1_{(p)}(z)$, that is a basis
of solutions to the corresponding qKZ equations.

Property 5 follows from the following property 6.

6.
$$
{\text{det}}_{1\leq m,\ell\leq n-1}\, (\, {2a_\ell\over p} \int_\R
G_m\Phi_pw_\ell dt)=
e^{(n-1)\pi i \sum^n_{j=1} z_j/p}
(2\pi i)^{n(n-1)/2} \times
$$
$$
\times
{\prod_{j=1}^n \Gamma({2a_j\over p} +1) \over \Gamma({2\over p}
\sum^n_{j=1}a_j \,+1)} \,\prod_{\ell<m}
\Gamma({z_m+a_m-z_\ell+a_\ell \over p})
\,\Gamma(1-{z_\ell+a_\ell-z_m+a_m \over p})
.
$$

{\it Example.} For $n=2$, we have Barnes' formula [WW]
$$
\int^{i\infty}_{-i\infty}
\Gamma(a+t)\Gamma(b+t)\Gamma(c-t)\Gamma(d-t)dt
=2\pi i{\Gamma(a+c)\Gamma(a+d)\Gamma(b+c)\Gamma(b+d)\over
\Gamma(a+b+c+d)}.
$$

\head 10. Conclusion\endhead

The theory of differential equations is a beautiful and
well developed theory. The
theory of difference equations is less elaborate.
One of the reasons is absence of well posed problems, good
examples, absence of indications to how to extend 
the notions of the theory of differential equations
to the theory of difference equations.
Mathematical physics (field theory, statistical 
mechanics) indicates such examples and a 
passage from differential to difference equations.
In particular, mathematical physics indicates a surprising quantization
of geometry in which algebraic functions are replaced by gamma functions,
such geometric objects  as cycles are replaced by
exponential functions and all basic relations
remain preserved.

\head References \endhead

\item{[F]} G. Felder, {\it Conformal field theory and integrable
systems associated to elliptic curves},
Proceedings of the International Congress of Mathematicians,
Zurich 1994, p. 1247--1255, Birkhauser, 1994;
 {\it Elliptic quantum groups,} preprint  hep-th/9412207,
to appear in the Proceedings of the ICMP, Paris 1994.

\item{[FTV]} G.Felder, V.Tarasov, A.Varchenko, {\it
Solutions to the elliptic qKZB equations and Bethe ansatz, I,}
(1996), preprint, 1-29.

\item{[FR]} I.Frenkel, N.Reshetikhin, {\it Quantum affine algebras
and holonomic difference equations,} Comm. Math. Phys. 146
(1992) 1-60.

\item{[JKMO]} M.Jimbo, A.Kuniba, T.Miwa, M.Okado, {\it
The $A^{(1)}_n$ face models,} Comm. Math. Phys., 119 (1988),543-565.

\item{[JM]} M.Jimbo, T.Miwa, {\it Algebraic analysis of solvable lattice
models,} AMS, 1995.

\item{[KRS]} P.Kulish, N.Reshetikhin, E.Sklyanin, {\it
Yang-Baxter equation and representation theory,} Lett. Math. Phys., 5
(1981), 393-403.

\item{[S]} F.Smirnov, {\it Form factors in completly integrable models of
quantum field theory,} World Scientific, 1992.

\item{[SV]} V.Schechtman, A.Varchenko, {\it Arrangements of
hyperplanes and Lie algebras homology,} Invent. Math. 106(1991),139-194.

\item{[TV1]} V.Tarasov, A.Varchenko, {\it
Asymptotic solutions to the Quantized Knizhnik-Zamolodchikov equation and
Bethe vectors,}  Amer. Math. Soc. Transl.(2), 174 (1996), 235-273.
                                                                       
\item{[TV2]} V.Tarasov, A.Varchenko, {\it
Geometry of $q$-Hypergeometric functions as a Bridge between
Yangians and Quantum Affine Algebras,}(1996), preprint, 1-84.

\item{[V]} A.Varchenko, {\it Multidimensional hypergeometric
functions and representation theory of Lie algebras and quantum groups,}
World Scientific, 1995.

\item {[WW]}  E.T.Whittaker and G.N.Watson,
{\it  A Course of Modern Analysis,}
Cambrige University Press , Cambribge, 1927.

\enddocument